\documentclass[12pt,a4paper]{amsart}
\usepackage{amsmath, amsfonts,amssymb, amsthm, amscd}
\input epsf.tex
\newcommand{\Z}{\mathbb Z}

\newtheorem{theo}{Theorem}
\newtheorem{cor}{Corollary}

\newtheorem{pro}{Proposition}
\title{On the singular braid monoid}
\author{V.~Vershinin}
\address{D\'epartement des Sciences Math\'ematiques,
Universit\'e Montpellier II,
Place Eug\'ene Bataillon,
34095 Montpellier cedex 5, France}
\email{ vershini@math.univ-montp2.fr}
\address{ Sobolev Institute of Mathematics, Novosibirsk, 630090,
Russia }
\email{ versh@math.nsc.ru}
\subjclass[2000]{Primary 20F36; Secondary 20F38}
\keywords{Braid group, singular braid monoid, word problem, Garside 
normal form. }
\thanks{
The author was supported
 in part by the French-Russian Program of
 Research EGIDE (dossier No 04495UL)}
\begin{document}
\begin{abstract} Garside's results and the existense of the greedy
normal form for braids are shown to be true
for the singular braid monoid. An analogue of the presentation
of J.~S.~Birman, K.~H.~Ko and S.~J.~Lee for the braid group is also
obtained for this monoid.
\end{abstract}
\maketitle
\section{Introduction}
Various questions concerning braid groups and their generalizations
attracted attention during the last decade. Presentations of the braid groups and algorithmic problems for the singular braid monoid
are among these questions

The canonical presentation of the braid group $Br_n$ was given by
E.~Artin \cite{Art1} and is well known. It has the generators
$\sigma_1$, $\sigma_2$, $\dots$, $\sigma_{n-1}$, and relations
\begin{equation}
 \begin{cases} \sigma_i \sigma_j &=\sigma_j \, \sigma_i, \ \
\text{if} \ \ |i-j|>1, \ \ i,j =1, ..., n-1;
\\ \sigma_i \sigma_{i+1} \sigma_i &= \sigma_{i+1} \sigma_i \sigma_{i+1},
\ \ i =1, ..., n-2.
\end{cases} \label{eq:brelations}
\end{equation}

Of course, there exist other presentations of the braid group.
J.~S.~Birman, K.~H.~Ko and S.~J.~Lee \cite{BKL} introduced a 
presentation with generators $a_{ts}$ with $1 \leq s<t\leq n$,
and relations

\begin{equation} \begin{cases}
a_{ts}a_{rq}&=a_{rq}a_{ts} \ \ {\rm for} \ \ (t-r)(t-q)(s-r)(s-q)>0,\\
a_{ts}a_{sr} &=a_{tr}a_{ts}=a_{sr}a_{tr}  \ \ {\rm for} \ \
1\leq r<s<t\leq n .
\end{cases}\label{eq:rebkl}
\end{equation}
The generators $a_{ts}$ are expressed by canonical generators
$\sigma_i$
in the following form:
 $$a_{ts}=(\sigma_{t-1}\sigma_{t-2}\cdots\sigma_{s+1})\sigma_s
(\sigma^{-1}_{s+1}\cdots\sigma^{-1}_{t-2}\sigma^{-1}_{t-1})  \ \
{\rm for} \ \ 1\leq s<t\leq n.$$

The {\it Baez--Birman monoid} $SB_n$  or {\it singular braid monoid}
\cite{Bae}, \cite{Bir2} is defined as a monoid with generators 
$\sigma_i, \sigma_i^{-1}, x_i$, $i=1,\dots,n-1,$ and relations
\begin{equation}
\begin{cases}
&\sigma_i\sigma_j=\sigma_j\sigma_i, \ \text {if} \ \ |i-j| >1,\\
&x_ix_j=x_jx_i, \ \text {if} \ \ |i-j| >1,\\
&x_i\sigma_j=\sigma_j x_i, \
\text {if} \ \ |i-j| \not=1,\\
&\sigma_i \sigma_{i+1} \sigma_i = \sigma_{i+1} \sigma_i
\sigma_{i+1},\\
&\sigma_i \sigma_{i+1} x_i = x_{i+1} \sigma_i \sigma_{i+1},\\
&\sigma_{i+1} \sigma_ix_{i+1} = x_i \sigma_{i+1} \sigma_i,\\
&\sigma_i\sigma_i^{-1}=\sigma_i^{-1}\sigma_i =1.
\end{cases}
\label{eq:singrel}
\end{equation}
In pictures $\sigma_i$ corresponds to the canonical generator of 
the braid group and $x_i$ represents an intersection
of the $i$th and $(i+1)$st strand as in
Figure~\ref{fi:singene}.
\begin{figure}
\begin{picture}(0,130)(0,-10) 
\thicklines
\put(0,50){\circle*{5}} \put(-100,100){\line(0,-1){100}}
\put(-50,100){\line(0,-1){100}} \put(-25,100){\line(1,-2){50}}
\put(25,100){\line(-1,-2){50}} \put(50,100){\line(0,-1){100}}
\put(100,100){\line(0,-1){100}}
\put(-100,110){\makebox(0,0)[cc]{$1$}}
\put(-50,110){\makebox(0,0)[cc]{$i-1$}}
\put(-25,110){\makebox(0,0)[cc]{$i$}}
\put(25,110){\makebox(0,0)[cc]{$i+1$}}
\put(50,110){\makebox(0,0)[cc]{$i+2$}}
\put(100,110){\makebox(0,0)[cc]{$n$}}
\put(-75,50){\makebox(0,0)[cc]{.\quad.\quad.}}
\put(75,50){\makebox(0,0)[cc]{.\quad.\quad.}}
\end{picture}
\caption{}\label{fi:singene}
\end{figure}
Motivation for introduction of this object was Vassiliev -- Goussarov
theory of finite type invariants.

The singular braid monoid on two strings is isomorphic to $\Z\oplus\Z^+$,
so the word problem in this case is trivial. For the
monoid with three strings this problem was solved by A.~J\'arai
\cite{Jar} and O.~T.~Dashbach and B.~Gemein \cite{DG}. In general
it was done by R.~Corran in
the complicated and technical paper \cite{Cor}. Here we hope to
give a simple solution using almost nothing but Garside's arguments.
The initial idea is simple and geometric: Garside's {\it fundamental
word} $\Delta$ (for convenience the definition is given below)
behaves the same way with respect to the singular generators $x_i$
as it behaves with respect to the braid generators $\sigma_i$:
\begin{equation*}
x_i\Delta \doteq \Delta x_{n-i},
\end{equation*}
as shown on Figure~\ref{fi:deltax}.
\begin{figure}
\begin{picture}(0,230)(0,-20)
\thicklines
\put(-120,200){\line(1,-1){12}} \put(-90,200){\line(-1,-1){30}}
\put(-60,200){\line(0,-1){30}} \put(-30,200){\line(0,-1){60}}
\put(50,200){\line(1,-1){30}} \put(80,200){\line(-1,-1){30}}
\put(110,200){\line(0,-1){60}} \put(140,200){\line(0,-1){90}}
\put(65,185){\circle*{5}}
\put(-104,182){\line(1,-1){24}} \put(-60,170){\line(-1,-1){30}}
\put(-120,170){\line(0,-1){60}}
\put(50,170){\line(1,-1){12}}
\put(80,170){\line(-1,-1){30}} \put(140,170){\line(0,-1){60}}
\put(66,152){\line(1,-1){24}}
\put(-74,152){\line(1,-1){24}}
\put(-90,140){\line(0,-1){30}}
\put(-30,140){\line(-1,-1){30}}
\put(50,140){\line(0,-1){60}} 
\put(110,140){\line(-1,-1){30}}
\put(96,122){\line(1,-1){24}}
\put(-90,110){\line(-1,-1){30}}
 \put(-120,110){\line(1,-1){12}}
\put(80,110){\line(0,-1){30}} 
\put(140,110){\line(-1,-1){30}}
\put(-30,110){\line(-1,1){12}} 
\put(-30,110){\line(0,-1){90}} 
\put(-102,92){\line(1,-1){24}} 
\put(-120,80){\line(0,-1){30}}
\put(-60,80){\line(-1,-1){60}} 
\put(-60,110){\line(0,-1){30}}
\put(50,80){\line(1,-1){12}} 
\put(80,80){\line(-1,-1){30}} 
\put(140,80){\line(-1,1){12}}
\put(140,80){\line(0,-1){90}}
\put(-72,62){\line(1,-1){12}} 
 \put(68,62){\line(1,-1){24}} 

\put(50,50){\line(0,-1){30}} 
\put(110,80){\line(0,-1){30}}
\put(-60,20){\line(1,-1){30}}
\put(50,20){\line(1,-1){12}} 
\put(-120,50){\line(1,-1){12}} 
\put(-90,20){\line(0,-1){30}} 
\put(-120,20){\line(0,-1){30}} 
\put(-30,20){\line(-1,-1){30}}
\put(-45,5){\circle*{5}}
\put(110,20){\line(0,-1){30}} 
\put(110,50){\line(-1,-1){60}}
\put(-60,50){\line(0,-1){30}} 
\put(98,32){\line(1,-1){12}}
\put(68,2){\line(1,-1){12}}
\put(-102,32){\line(1,-1){12}}
\put(10,90){\makebox(0,0)[cc]{$-$}}
\put(10,85){\makebox(0,0)[cc]{$-$}}
\end{picture}
\caption{}\label{fi:deltax}
\end{figure}

\section{Word problem for the singular braid monoid}

 Following Garside's ideas
we consider the {\it positive singular braid monoid} $SB_n^+$. It is
defined as a monoid with generators $\sigma_i, x_i$, $i=1,\dots,n-1,$ and
relations (\ref{eq:singrel}) except the last one concerning the
invertibility of $\sigma_i$. Two positive words $A$ and $B$ in the
alphabet $\{\sigma_i, x_i$, $(i=1,\dots,n-1) \}$
will be said to be {\it positively equal} if they are equal as elements
of $SB_n^+$. In this case we shall write $A\doteq B$. As usual, identity
of words is denoted by the symbol $\equiv$.
Proofs of statements below are the same as in Garside's paper
\cite{Gar} with some exceptions as, for example, the proof of
Proposition~\ref{pro:inj}. Garside's proof of this Proposition
doesn't work and we use the
Malcev rule. Proposition~\ref{pro:inj} was proved by R.~Corran \cite{Cor},
Theorem~\ref{theo:center} was proved by
R.~Fenn, D.~Rolfsen and J.~Zhu \cite{FRZ} using different methods.
\begin{pro} For $i,k = 1, ..., n-1 $, given
$\sigma_i A \doteq \sigma_k B$,
it follows that
\begin{equation*} \text{if} \ \  k = i, \ \text{then} \ \ A \doteq B,
\end{equation*}
\begin{equation*} \text{if} \ \  |k - i| = 1,
\ \text{then} \ \ A \doteq \sigma_k\sigma_i Z, \ B\doteq \sigma_i\sigma_k Z
\ \ \text{for some} \ \ Z,
\end{equation*}
\begin{equation*} \text{if} \ \ |k - i| \geq 2, \ \text{then} \ \ A \doteq
\sigma_k Z, \  B\doteq \sigma_i Z \ \ \text{for some} \ \ Z,
\end{equation*}
given
$\sigma_i A \doteq x_k B$,
it follows that
\begin{equation*} \text{if} \ \  |k - i| = 1, 
\ \text{then} \ \ A \doteq \sigma_k x_i Z, \ B\doteq \sigma_i\sigma_k Z
\ \ \text{for some} \ \ Z,
\end{equation*}
\begin{equation*} \text{if} \ \ |k - i| \not= 1, \ \text{then} \ \ A 
\doteq
x_k Z, \ B\doteq \sigma_i Z \ \ \text{for some} \ \ Z.
\end{equation*}
and given
$x_i A \doteq x_k B$,
it follows that
\begin{equation*} \text{if} \ \  k = i, \ \text{then} \ \ A \doteq B,
\end{equation*}
\begin{equation*} \text{if} \ \  |k - i| \geq 2, 
\ \text{then} \ \ A \doteq x_k Z, \ B\doteq x_i Z
\ \ \text{for some} \ \ Z,
\end{equation*}
\begin{equation*} \text{the case when } \  |k - i| = 1, \ 
\text{is impossible}.
\end{equation*}
The same is true for multiples of $\sigma_i$ or $x_k$ to the right.
\label{pro:div}
\end{pro}
\begin{proof}
Garside's proof works here. We apply the induction on the length $s$ 
of $A$ and the length of chain of transformations from $a_iA$ to $a_kB$
where $a_i$ may be $\sigma_i$ or $x_i$ and $a_k$ may be $\sigma_k$
or $x_k$.
The cases of $s=0$, $1$ are evident, so suppose that the statement is 
true for length $s\leq r$ and for $s=r+1$ it is true for chain-length
$\leq t$. As an example we give a proof of the last 
statement, which is formally is not contained in Garside's considerations.
So, let $A$, $B$ be of word-length $r+1$ and let $x_i A \doteq x_kB$, 
$|i-k|=1$, through a transformation of chain-length $t+1$. We may
suppose that $k=i+1$ and let the succesive words of transformations be
\begin{equation*}
W_1 \equiv x_iA, \dots W_{t+2}\equiv x_{i+1}B.
\end{equation*}
Choose arbitrary any intermediate word $W_g$, say, from the middle of
the chain somewhere. We have $W\equiv aV$, where $a$ is a generator of 
$SB_n^+$. Suppose at first that $a$ commutes with $x_i$, $x_{i+1}$, 
then we have $x_i A\doteq a V\doteq x_{i+1}B$ and using induction we
obtain:
\begin{equation*}
A\doteq a P, \ \ V\doteq x_i P, \ \ V\doteq x_{i+1} Q, \ \ B\doteq aQ.
\end{equation*}
So, $x_i P \doteq x_{i+1} Q$ what is impossible by induction.
 
The cases when $a= x_{i-1}, x_i, x_{i+1}, x_{i+2} $ are also impossible
by induction. The cases which we need to consider are $a= \sigma_{i-1}, 
\sigma_i, \sigma_{i+1}, \sigma_{i+2} $. So, let $a=\sigma_{i-1}$, then
$W_g\doteq \sigma_{i-1}V$. Using induction we get
\begin{equation*}
A\doteq \sigma_{i-1}\sigma_i P, \ \ V\doteq \sigma_i x_{i-1} P, \ \ 
V\doteq x_{i+1} Q, \ \ B\doteq \sigma_{i-1}Q.
\end{equation*}
Hence $\sigma_i x_{i-1} P \doteq x_{i+1} Q$, and so, $ x_{i-1} P \doteq \sigma_{i+1}x_i R$, $Q\doteq\sigma_i\sigma_{i+1}R$. Again using induction 
we get $P\doteq \sigma_{i+1}S$, $x_iR\doteq x_{i-1}S$, that is impossible.
The rest cases may be considered the same way. 
\end{proof}
\begin{cor} If $A\doteq P$, $B\doteq Q$, $AXB \doteq PYQ$, ($L(A)\geq 0$,
$L(B) \geq 0$), then $X\doteq Y$. That is, monoid $SB_n$ is left and right
cancellative.
\end{cor}
\begin{pro} For any word $W$ in the alphabet $\{\sigma_i$, $x_i$ 
$(i = 1, ..., n-1 ) \}$, let $S$ be a
word in the alphabet $\{\sigma_i\}$ of maximal length, such that $W\doteq ST$
for some word $T$. Suppose also that $W\doteq AV$ for some word $A$ in the
alphabet $\{\sigma_i\}$. Then $S$ is divisible by $A$. The same is true for 
the right division.
\label{pro:maxdiv}
\end{pro}
\begin{proof}We use the induction on the length of $S$. If $S$ has the 
length 1, then $A$ also has the length 1 and the assertion follows from 
Proposition~\ref{pro:div}. Let the statement be true for the length less or 
equal to $k$ and
let the length of $S$ be equal to $k+1$. 
Consider at first the case when the length of $A$ is equal to 1. This means 
that $A\equiv \sigma_j$. Let the first letter of $S$ be $\sigma_i$: 
$S\equiv \sigma_i S^\prime$, so we have 
the situation $\sigma_i S^\prime T\doteq \sigma_j V$. If $i=j$, then 
we are done.
If $|i-j| \geq 2$, then $S^\prime T \doteq \sigma_j X$ and using the induction
we have $S^\prime  \doteq \sigma_j R^\prime$. Let $|i-j|=1$, then   
$S^\prime T \doteq \sigma_j\sigma_i Y$ and again using induction we obtain
$S^\prime  \doteq \sigma_j\sigma_i Q^\prime$. Suppose now that the 
length of $A$ 
is greater than l, then $A\equiv \sigma_j A^\prime$ for some $\sigma_j$.
Using the previous case
we have $S\doteq \sigma_j S^\prime$ and $S^\prime T \doteq A^\prime V$ 
and by
induction we obtain $S^\prime\doteq A^\prime S^{\prime\prime}$. Hence $S$ is 
divisible by $A$.
\end{proof}
\begin{cor} For any element $w$ of $SB_n^+$ there exists a unique greatest left
(right) diviser which belongs to $Br_n^+$. 
\end{cor}
Garside's {\it fundamental word} for the braid group $Br_n$ is
the following
$$\Delta \equiv \sigma_1 \dots \sigma_{n-1}\sigma_1 \dots \sigma_{n-2}
\dots \sigma_1\sigma_2\sigma_1.$$
If we use Garside's notation $\Pi_t\equiv \sigma_1\dots \sigma_t$, then
$\Delta \equiv \Pi_{n-1} \dots \Pi_1$.
We keep the same notations for their images in $SB_n$.
Garside's transformation of words $\mathcal R$ and then the
automorphism of $Br_n$ and the positive braid monoid $Br_n^+$,
defined by the formula
\begin{equation*}
\mathcal R(\sigma_i) \equiv \sigma_{n-i},
\end{equation*}
we extend to letters $x_i$ and so to $SB_n^+$ and to $SB_n$ by
\begin{equation*}
\mathcal R(x_i) \equiv x_{n-i}.
\end{equation*}
\begin{pro} There are equalities
$$\sigma_i\Delta \doteq\Delta \mathcal R(\sigma_i),$$
$$x_i\Delta \doteq\Delta \mathcal R(x_i).$$
\end{pro}
\begin{proof}
\begin{multline*}
x_1\Delta \equiv x_1\sigma_1 \dots \sigma_{n-1}\sigma_1 \dots \sigma_{n-2}
\dots \sigma_1\sigma_2\sigma_1
\doteq \sigma_1 x_1\sigma_2 \dots \sigma_{n-1}\sigma_1 \dots \sigma_{n-2}
\dots \sigma_1\sigma_2\sigma_1
\doteq \\
\sigma_1 x_1\sigma_2 \sigma_1\dots \sigma_{n-1}\sigma_2 \dots \sigma_{n-2}
\dots \sigma_1\sigma_2\sigma_1
\doteq \sigma_1 \sigma_2 \sigma_1 x_2\sigma_3\dots \sigma_{n-1}\sigma_2 \dots \sigma_{n-2}
\dots \sigma_1\sigma_2\sigma_1
\doteq \\
\sigma_1 \sigma_2 \sigma_1 \sigma_3\dots x_{n-2}\sigma_{n-1}\sigma_{n-2}
\dots \sigma_1\sigma_2\sigma_1
\doteq \sigma_1 \sigma_2 \sigma_1 \sigma_3\dots \sigma_{n-2}\sigma_{n-1}x_{n-1}
\dots \sigma_1\sigma_2\sigma_1
\doteq \\
\sigma_1 \sigma_2  \sigma_3\dots \sigma_{n-2}\sigma_{n-1}\sigma_1 \dots\sigma_{n-2}
\dots \sigma_1\sigma_2\sigma_1 x_{n-1}\equiv \Delta x_{n-1}.
\end{multline*}
Let $i\geq 2$:
\begin{multline*}
x_i\Delta \equiv x_i\sigma_1 \dots \sigma_{n-1}\sigma_1 \dots \sigma_{n-2}
\dots \sigma_1\sigma_2\sigma_1
\doteq \sigma_1 \dots x_i\sigma_{i-1}\sigma_i \dots \sigma_{n-1}\sigma_1 \dots \sigma_{n-2}
\dots \sigma_1\sigma_2\sigma_1
\doteq \\
\sigma_1 \dots \sigma_{i-1}\sigma_i x_{i-1}\dots \sigma_{n-1}\sigma_1 \dots \sigma_{n-2}
\dots \sigma_1\sigma_2\sigma_1
\doteq
\Pi_{n-1}x_{i-1}\Pi_{n-2} \dots\Pi_1 \doteq \\
\Pi_{n-1}\dots\Pi_{n-i+1}x_{1}\Pi_{n-i}\Pi_{n-i-1} \dots\Pi_1 \doteq
\Pi_{n-1}\dots\Pi_{n-i+1}\Pi_{n-i}\Pi_{n-i-1}x_{n-i}\Pi_{n-i-2} \dots\Pi_1 \doteq \\
\Pi_{n-1}\dots\Pi_{n-i+1}\Pi_{n-i} \dots\Pi_1 x_{n-i}\doteq
\Delta x_{n-i}.
\end{multline*}
\end{proof}
\begin{pro} If $W$ is any positive word in $SB_n^+$ such that either
\begin{equation*}
W\doteq \sigma_1 A_1\doteq\sigma_2 A_2\doteq\dots\doteq\sigma_{n-1}A_{n-1},
\end{equation*}
or
\begin{equation*}
W\doteq B_1\sigma_1 \doteq B_2\sigma_2 \doteq\dots\doteq
B_{n-1}\sigma_{n-1},
\end{equation*}
then $W\doteq\Delta Z$ for some $Z$.
\end{pro}
\begin{pro} Canonical homomorphism
$$SB_n^+ \to SB_n$$
is a monomorphism.
\label{pro:inj}
\end{pro}
\begin{proof} We need to prove that if two elements of $SB_n^+$,
(expressed by positive words $A$ and $C$) are equal in $SB_n$ then
they are positively equal. If elements defined by words $A$ and $B$
are equal in $SB_n$ then there exists a sequence of words
\begin{equation}
A\equiv A_0 \to A_1 \to \dots \to A_j \to\dots \to A_k\equiv C,
\label{eq:seq}
\end{equation}
where each arrow means an elementary operation which may be an
application of one defining relation or insert or deletion of an
expression $c c^{-1}$ or $c^{-1}c$ where $c$ is one of $\sigma_i$.
In the first case the insert is called \emph{left} and in the second
it is called \emph{right}. A right insert
\begin{equation*}
A_j \equiv Y_jZ_j \to  Y_j c^{-1}cZ_j
\end{equation*}
is called \emph{correct}, if the fragment $Y_j$ is not changed in the
sequence (\ref{eq:seq}) till the elimination of $c^{-1}$. For the left
insert the same condition is made for the fragment $Z_j$.

We use the following Malcev rule \cite {Mal1}, \cite{Mal2},
\cite{Nov}: 

\emph{If all defining relations contain only positive degrees of
letters of the alphabet then transformation of a word $A$ into
a positive word $C$ can be done by a sequence of operations where all
inserts are correct.}

Let us consider a subsequence of (\ref{eq:seq}) starting with
the last right insert and finishing with the deletion of the
corresponding $c^{-1}$:
\begin{equation*}
A_j\equiv Y_jZ_j \to  Y_j c^{-1}cZ_j \to Y_j c^{-1}V_{j+1}\to
\dots \to Y_j c^{-1}cZ_{j+r}\to Y_j Z_{j+r}.
\end{equation*}
The fragment $Y_j$ does not change in this subsequence because the
insert is correct. For any $c$ there exists a positive word $D_c$
such that $D_c c\doteq \Delta$. Define the sequence
\begin{multline*}
\Delta A_j\equiv \Delta Y_jZ_j \mapsto \dots \mapsto
\mathcal R (Y_j) \Delta Z_j \mapsto \dots \mapsto
\mathcal{R}(Y_j) D_c c Z_{j} \hookrightarrow
\mathcal{R}(Y_j) D_c V_{j+1} \hookrightarrow \dots \\
\hookrightarrow \mathcal{R}(Y_j) D_c c Z_{j+r}
\mapsto \dots \mapsto \mathcal{R}(Y_j) \Delta Z_{j+r}\mapsto
\dots \mapsto \Delta Y_j Z_{j+r}
\end{multline*}
where $\mapsto$ denote a positive operation and $\hookrightarrow$
denote a positive operation or deletion. So multiplying by $\Delta$
we eliminated one insert. Applying induction we construct a
sequence of positive operations between $\Delta^m A\Delta^l$ and
$\Delta^m C\Delta^l$ for some $m$ and $l$. After the cancellation this
means that $A$ and $C$ are positive equivalent. As it was noted 
by V.~V.~Chainikov,
one can get a proof without use of the Malcev rule but using the fact that
relations (\ref{eq:singrel}) are invariant with respect to the operation
$\mathcal{R}$.
\end{proof}
Among all positive words on the alphabet  $\{\sigma_1, \dots, \sigma_{n-1}$,
$x_1, \dots, x_{n-1} \}$ let us
introduce a lexicographical ordering with the condition that
$\sigma_1 < \sigma_2 < \dots < \sigma_{n-1} < x_1 < x_2 < \dots <x_{n-1} $.
For a positive word $W$ the
\emph{base} of $W$ is the smallest positive word with respect to this
ordering, which is positively equal
to $W$. The base is uniquely determined. If a positive word $A$ is prime
to $\Delta$, then for the base of $A$ the notation $\overline{A}$  will
be used.
\begin{theo} In $SB_n$ every word $W$ can be expressed uniquely in the
form $ \Delta^m \overline{A}$, where $m$ is an  integer.
\label{theo:garnf}
\end{theo}
\begin{proof} First suppose $P$ is any positive word. Among all positive
words positively equivalent to $P$ choose a word in the form
$\Delta^t A$ with $t$ maximal. Then $A$ is prime to $\Delta$ and we have
\begin{equation*}
P\doteq \Delta^t \overline{A}
\end{equation*}
Now let $W$ be any word in $SB_n$. Then we may put
\begin{equation*}
W\equiv W_1(c_1)^{-1}W_2(c_2)^{-1} \dots (c_k)^{-1}W_{k+1},
\end{equation*}
where each $W_j$ is a positive word of length $\geq 0$, and $c_l$
are generators $\sigma_i$, the only possible invertable generators. 
As it was already mentioned for each $c_l$ there exists a positive 
word $D_l$ such that $c_l D_l\doteq \Delta$, so that
$(c_l)^{-1} = D_l \Delta^{-1}$, and hence
\begin{equation*}
W = W_1 D_1 \Delta^{-1}W_2 D_2 \Delta^{-1} \dots
W_k D_k \Delta^{-1}W_{k+1}.
\end{equation*}
Hence, moving the factors $\Delta^{-1}$ to the left, we obtain
$W=\Delta^k P$, where $P$ is positive, so we can express it in the
form $\Delta^t\overline{A}$ and finally we get
\begin{equation}
W = \Delta^m \overline{A}.
\label{eq:gnf}
\end{equation}
It remains to show that the form (\ref{eq:gnf}) is unique. Suppose
\begin{equation}
 \Delta^m \overline{A}= \Delta^p \overline{C}.
\label{eq:gnf2}
\end{equation}
Let $p<m$, and $m-p=t>0$. Then (\ref{eq:gnf2}) gives
$\Delta^t\overline{A}= \overline{C}$, what is impossible. So $p =m$
and hence $\overline{A} = \overline{B}$. So from
Proposition~\ref{pro:inj} we obtain $\overline{A}\doteq\overline{C}$,
but the base is unique, hence $\overline{A}\equiv\overline{C}$ and
the uniqueness of the form (\ref{eq:gnf}) is established.
\end{proof}
The form of a word $W$ established in this theorem we call the
\emph{Garside left normal form} and the 
index $m$ we call the \emph{power} of $W$. The same way the
\emph{Garside right normal form} is defined and the corresponding
variant of Theorem~\ref{theo:garnf} is true. The
Garside normal form also gives a solution to the word problem in the braid
group.
\begin{cor}
The necessary and sufficient condition that two words in $SB_{n}$ 
are equal is that their Garside normal forms (left or right) are identical.
\label{cor:gws}
\end{cor}
Garside normal form for the braid groups was precised in the subsequent 
papers \cite{Ad}, \cite{E_Th}, \cite{EM}. Namely, it
was written in the \emph{left-greedy form} (in the terminology of
W.~Thurston \cite{E_Th}) 
\begin{equation*}
\Delta^t A_1 \dots A_k,
\end{equation*}
where $A_i$ are the succesive
possible longest \emph{fragments of the word} $\Delta$ (in the terminology 
of S.~I.~Adyan \cite{Ad}) or \emph{positive permutation braids} (in the 
terminology of E.~El-Rifai and H.~R.~Morton \cite{EM}). Certainly, the same 
way the \emph{right-greedy form} is defined. 

Consider these forms for the singular braid monoid.  
For any word $W$ first of all move to the left the greatest power of $\Delta$.
To the right we have a positive word $W^\prime$ not divisible by $\Delta$.
Consider the 
decomposition $W^\prime\doteq S_1T$ of Proposition~\ref{pro:maxdiv}. 
Then we take
the fragments of the left-greedy form for 
$S_1$: $S_{1,1} \dots S_{1,t}$. Among all the 
$x_i$-divisers of $T$ we choose the smallest in the lexicographical order:
$T\doteq x_{i_1} T_1$. Consider the decomposition of 
Proposition~\ref{pro:maxdiv} for $T_1$ and continue this process. We get 
the form
\begin{equation*}
W\doteq \Delta^t S_1 X_1 \dots S_kX_k,
\end{equation*}
where each $S_i$ consists of the fragments of the left-greedy form for the braid
group and each $X_i$ is a 
lexicographically ordered product of $x_j$. We call this form the 
\emph{left-greedy form} for the singular braid monoid. The same way
the \emph{right-greedy form} for the singular braid monoid is defined.
\begin{theo} Each element of the singular braid monoid can be written uniquely
in a left-greedy (right-greedy) form. 
\end{theo}

\begin{theo} For $n=2$ the singular braid monoid $SB_n$ is commutative
and isomorphic to $\Z\oplus\Z^+$. For $n\geq 3$ its center is the same
as the center of $Br_n$ and so is generated by $\Delta^2$.
\label{theo:center}
\end{theo}

\section{Conjugacy problem for the singular braid monoid}

Let $M$ be a monoid with unit group $G$. We call two elements $u, v\in M$
\emph{conjugate} if $v = g^{-1} u g$ for some $g\in G$. This means that
the elements $u$ and $v$ belong to the same orbit of the canonical
action of the braid group $Br_n$ on $SB_n$. So, in this sense
we understand the \emph{conjugacy problem} for the singular braid monoid.
Such approach to conjugacy appears in monoid theory, see
\cite{Put}, for example.

Relations (\ref{eq:singrel}) are homogeneous with respect to both types of
generators $\sigma_i$ and  $x_i$, so, three homomorphisms which express the
degree of an element with respect to $\sigma_i$, $x_i$ and the total
degree are well defined:
$\operatorname{deg}_\sigma :SB_n\to \Z$,
$\operatorname{deg}_x :SB_n\to \Z^+$ and
$\operatorname{deg} :SB_n\to \Z$.

\begin{pro} The group of units of monoid $SB_n$ is equal to the image of
the braid group $Br_n$.
\end{pro}
\begin{proof} Invertible elements must have invertible degree
$\operatorname{deg}_x$,
so $\operatorname{deg}_x = 0$.
\end{proof}
Garside's solution of the conjugacy problem for the braid groups works
in the case of the singular braid monoid. In the proof
of the following theorem concerning the structure of the Cayley diagram
$D(W)$ of a positive word $W$ on the alphabet $\{\sigma_i, x_k\}$,
containing $\Delta$, some additional considerations are necessary.

\begin{theo} If $W \doteq \Delta V$ is any positive word (on the alphabet
$\sigma_i$, $x_i$) containing $\Delta$, then each node of $D(W)$ is
incident with each edge $\sigma_i$, $i=1,\dots n-1$.
\end{theo}
\begin{proof} To complete Garside's proof by induction on the order of
a node we must consider the case when
the edge $x_i$ starts at some node $C$ of order $m$ and ends at a node
$D$ of order $m+1$.

We first consider the case of the generators $\sigma_k$, with
$|k -i| \not= 1$. If the edge $\sigma_k$ ends at $C$ then we have
$\sigma_k x_i \doteq x_i\sigma_k$ what means that the edge $\sigma_k$
also ends at $D$. If $\sigma_k$ starts from $C$ then using
Proposition~\ref{pro:div} we obtain the fragment of the Cayley graph,
depicted on Figure~\ref{fi:x_sigabc} a).
\begin{figure}
\epsfbox{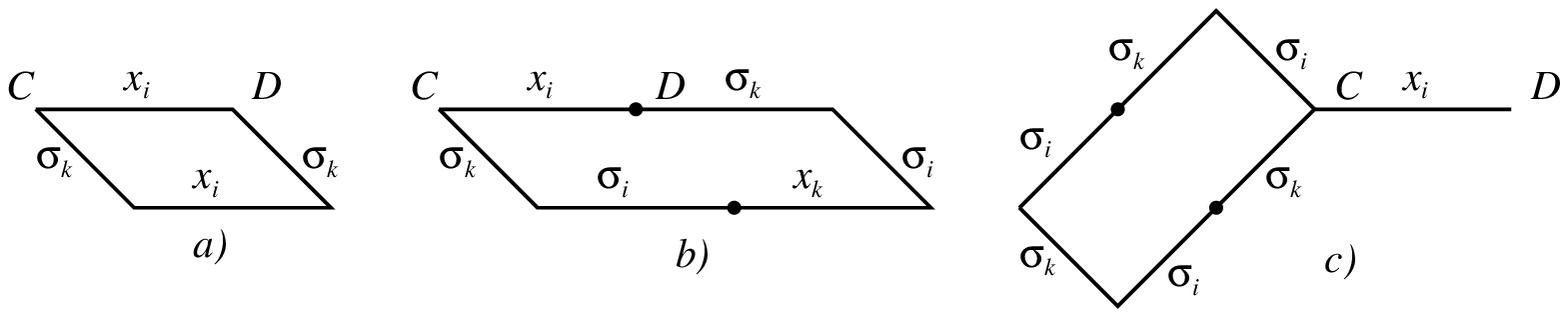}
\caption{}\label{fi:x_sigabc}
\end{figure}
\begin{figure}
\epsfbox{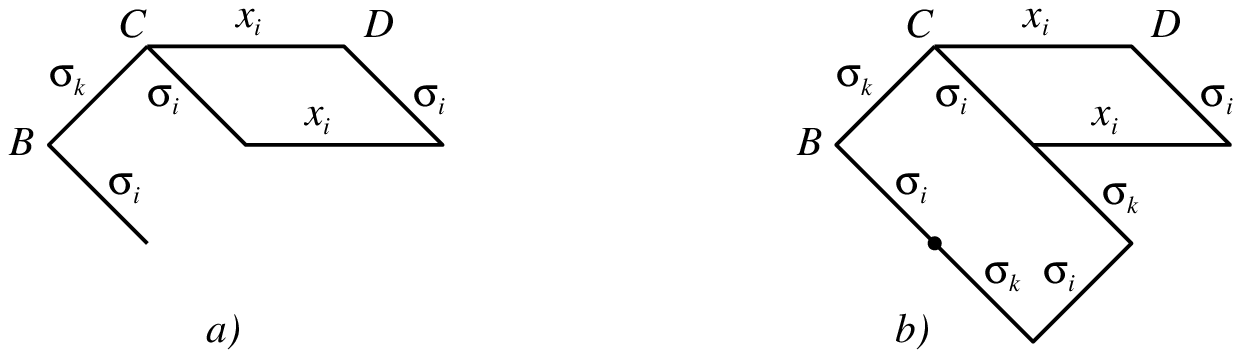}
\caption{}\label{fi:x_sidab}
\end{figure}
\begin{figure}
\epsfbox{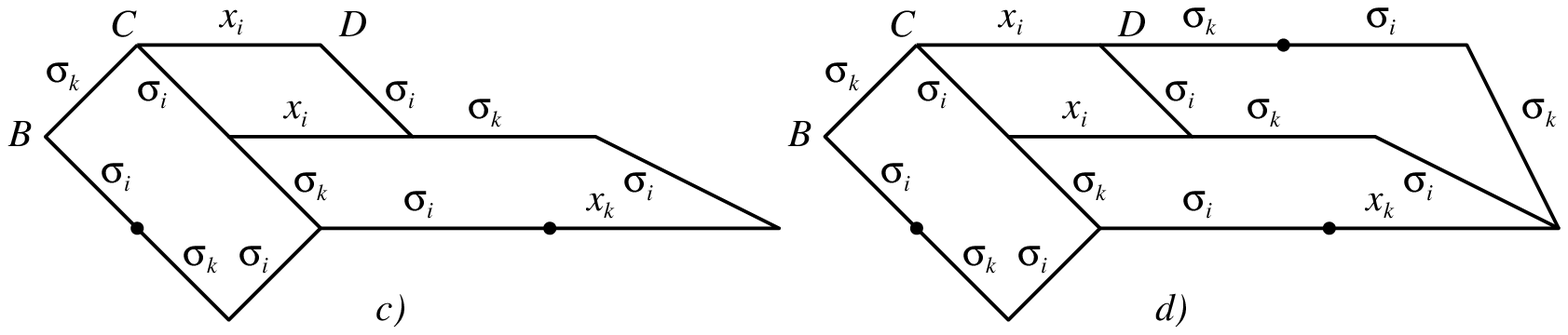}
\caption{}\label{fi:x_sidcd}
\end{figure}

Let $|k-i| =1$. If the edge $\sigma_k$ starts at $C$, then using
Proposition~\ref{pro:div} we obtain the fragment of the Cayley
graph, depicted on Figure~\ref{fi:x_sigabc} b).
Consider now the case when the edge $\sigma_k$ ends at $C$ (and starts
at some node $B$). The edge
$\sigma_i$ must be incident to $C$, so it either ends or starts at $C$.
Consider the first case, then using Proposition~\ref{pro:div}
we obtain the fragment of the Cayley graph,
depicted on Figure~\ref{fi:x_sigabc} c).
We complete the proof of this case using relation
$\sigma_i\sigma_k x_i= x_{k}\sigma_i\sigma_k$.
At the last case the edge $\sigma_i$ starts at $C$.
Consider the node $B$ where starts $\sigma_k$. If the edge $\sigma_i$
ends at $B$ then we  use relation
$\sigma_i\sigma_k x_i= x_{k}\sigma_i\sigma_k$ and this case
is finished. So suppose that the edge $\sigma_i$ also starts at $B$.
Using
Proposition~\ref{pro:div} we obtain the fragment of the Cayley graph,
depicted on Figure~\ref{fi:x_sidab} a). Using
Proposition~\ref{pro:div} two times we obtain the fragment of the
Cayley graph, depicted on Figure~\ref{fi:x_sidab} b).
Again using
Proposition~\ref{pro:div} we come to the fragment of the Cayley graph,
depicted on Figure~\ref{fi:x_sidcd} c).
To complete the proof of this case we use the relation
$\sigma_i\sigma_k x_i= x_{k}\sigma_i\sigma_k$ and arrive to
Figure~\ref{fi:x_sidcd} d).
\end{proof}
The definition of the summit set for the singular braid monoid is the
same as Garside's.
\begin{theo} In $SB_n$ two elements are conjugate if and only if their
summit sets are identical.
\end{theo}
For any element of $SB_n$ the sumit set is finite and is obtained 
algorithmically, so this theorem gives a solution of the conjugacy 
problem.

\section{Birman -- Ko -- Lee presentation for the singular
braid monoid}

\begin{figure}
\epsfbox{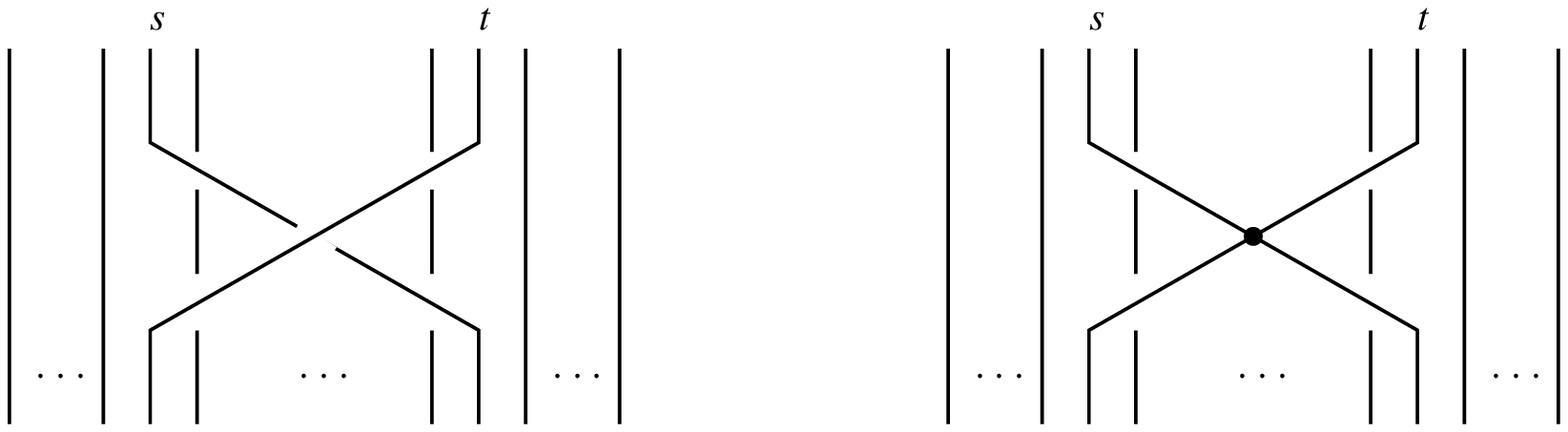}
\caption{}
\label{fi:sbige}
\end{figure}

For the singular braid monoid we prove the existence of the analogue
of the presentation of
J.~S.~Birman,  K.~H.~Ko and S.~J.~Lee (\ref{eq:rebkl}).
For $1\leq s<t\leq n$ and
$1\leq p<q\leq n$ we consider the elements of $SB_n$ which are
defined by
\begin{equation} \begin{cases}
a_{ts}&=(\sigma_{t-1}\sigma_{t-2}\cdots\sigma_{s+1})\sigma_s
(\sigma^{-1}_{s+1}\cdots\sigma^{-1}_{t-2}\sigma^{-1}_{t-1})  \ \
{\rm for} \ \ 1\leq s<t\leq n, \\
a_{ts}^{-1}&=(\sigma_{t-1}\sigma_{t-2}\cdots\sigma_{s+1})\sigma_s^{-1}
(\sigma^{-1}_{s+1}\cdots\sigma^{-1}_{t-2}\sigma^{-1}_{t-1})  \ \
{\rm for} \ \ 1\leq s<t\leq n, \\
b_{qp}&=(\sigma_{q-1}\sigma_{q-2}\cdots\sigma_{p+1}) x_p
(\sigma^{-1}_{p+1}\cdots\sigma^{-1}_{q-2}\sigma^{-1}_{q-1})  \ \
{\rm for} \ \ 1\leq p<q\leq n.
\end{cases}\label{eq:defab}
\end{equation}
Geometrically the generators $a_{s,t}$ and $b_{s,t}$ are depicted
in Figure~\ref{fi:sbige}.
\begin{theo} The singular braid monoid $SB_n$ has a presentation with
generators $a_{ts}$, $a_{ts}^{-1}$
for $1\leq s<t\leq n$ and  $b_{qp}$ for
$1\leq p<q\leq n$
and relations
\begin{equation} \begin{cases}
a_{ts}a_{rq}&=a_{rq}a_{ts} \ \ {\rm for} \ \ (t-r)(t-q)(s-r)(s-q)>0,\\
a_{ts}a_{sr} &=a_{tr}a_{ts}=a_{sr}a_{tr}  \ \ {\rm for} \ \
1\leq r<s<t\leq n , \\
a_{ts}a_{ts}^{-1} &=a_{ts}^{-1}a_{ts} =1 \ \ {\rm for} \ \ 1\leq s<t\leq n,\\
a_{ts}b_{rq}&=b_{rq}a_{ts} \ \ {\rm for} \ \ (t-r)(t-q)(s-r)(s-q)>0,\\
a_{ts}b_{ts}&=b_{ts}a_{ts}  \ \ {\rm for} \ \
1\leq s<t\leq n , \\
a_{ts}b_{sr} &=b_{tr}a_{ts}  \ \ {\rm for} \ \
1\leq r<s<t\leq n , \\
a_{sr}b_{tr} &=b_{ts}a_{sr}  \ \ {\rm for} \ \
1\leq r<s<t\leq n , \\
a_{tr}b_{ts}&=b_{sr}a_{tr}  \ \ {\rm for} \ \
1\leq r<s<t\leq n, \\
b_{ts}b_{rq}&=b_{rq}b_{ts} \ \ {\rm for} \ \ (t-r)(t-q)(s-r)(s-q)>0.
\end{cases}\label{eq:srebkl}
\end{equation}
\end{theo}
\begin{proof} We follow 
the proof of
J.~S.~Birman, K.~H.~Ko and S.~J.~Lee \cite{BKL}
and we begin with the presentation of $SB_n$ using generators
$\sigma_i, \sigma_i^{-1}, x_i$, $i=1,\dots,n-1,$ and relations
(\ref{eq:singrel}). Add the new generators
$a_{ts}$, $a_{ts}^{-1}$
for $1\leq s<t\leq n$ and  $b_{qp}$ for
$1\leq p<q\leq n$ and relations (\ref{eq:defab}).
Relations (\ref{eq:srebkl}) are described by isotopies of
singular braids,
so they must be the consequences of (\ref{eq:singrel}), and we may
add them too.
\par
In the special case when $t = s+1$ relations (\ref{eq:defab}) tell
us that
$
a_{(s+1) s} =\sigma_s,
$
$
a_{(s+1) s}^{-1} =\sigma_s^{-1},
$
$
b_{(s+1) s} =x_s,
$
so we may omit the generators $\sigma_1, \dots, \sigma_{n-1}$,
$\sigma_1^{-1}, \dots, \sigma_{n-1}^{-1}$
to obtain a presentation with generators $a_{ts}$, $a_{ts}^{-1}$,
$b_{pq}$. Defining relations are now
(\ref{eq:srebkl}) and
\begin{equation}
a_{(i+1)i} a_{(j+1)j}= a_{(j+1)j} a_{(i+1)i}, \ \text {if} \ \ |i-j| >1,
\label{eq:comma}\end{equation}
\begin{equation}
a_{(i+1)i} a_{(i+2)(i+1)} a_{(i+1)i} = a_{(i+2)(i+1)} a_{(i+1)i}
a_{(i+2)(i+1)},
\label{eq:brel}
\end{equation}
\begin{equation}
a_{ts}=(a_{t(t-1)}a_{(t-1)(t-2)}\cdots a_{(s+2)(s+1)})a_{(s+1)s}
(a^{-1}_{(s+2)(s+1)}\cdots a^{-1}_{(t-1)(t-2)}a^{-1}_{t(t-1)})
\label{eq:defaa} \end{equation}
\begin{equation}
a_{ts}^{-1}=(a_{t(t-1)}a_{(t-1)(t-2)}\cdots a_{(s+2)(s+1)})
a_{(s+1)s}^{-1}
(a^{-1}_{(s+2)(s+1)}\cdots a^{-1}_{(t-1)(t-2)}a^{-1}_{t(t-1)})
\label{eq:defaa-1} \end{equation}
\begin{equation}
b_{(i+1)i}b_{(j+1)j}=b_{(j+1)}b_{(i+1)i}, \ \text {if} \ \ |i-j| >1,
\label{eq:commb}\end{equation}
\begin{equation}
b_{(i+1)i} a_{(j+1)j}= a_{(j+1)j} b_{(i+1)i}, \
\text {if} \ \ |i-j| \not=1,
\label{eq:commab}\end{equation}
\begin{equation}
a_{(i+1)i} a_{(i+2)(i+1)} b_{(i+1)i} = b_{(i+2)(i+1)} a_{(i+1)i}
a_{(i+2)(i+1)},
\label{eq:brelab1}\end{equation}
\begin{equation}
a_{(i+2)(i+1)} a_{(i+1)i} b_{(i+2)(i+1)} = b_{(i+1)i} a_{(i+2)(i+1)}
a_{(i+1)i},
\label{eq:brelab2}\end{equation}
\begin{multline}
b_{qp}=(a_{q(q-1)}a_{(q-1)(q-2)}\cdots a_{(p+2)(p+1)}) b_{(p+1)p}
(a^{-1}_{(p+2)(p+1)}\cdots a^{-1}_{(q-1)(q-2)} a^{-1}_{q(q-1)}) \\
 {\rm for} \ \ 1\leq p<q\leq n.
\label{eq:defbab}
\end{multline}
Now we prove that relations (\ref{eq:comma}) - (\ref{eq:defbab})
are consequences of (\ref{eq:srebkl}).
It was proved by J.~S.~Birman, K.~H.~Ko and S.~J.~Lee that
relations (\ref{eq:comma}) - (\ref{eq:defaa}) are consequences
of the first two relations of (\ref{eq:srebkl}). The proof for
the relation (\ref{eq:defaa-1}) is the same as for (\ref{eq:defaa}).
Relations (\ref{eq:commb}) are special cases of the last relations
in (\ref{eq:srebkl}).
Relations (\ref{eq:commab}) are special cases of the forth and
the fifth relations in (\ref{eq:srebkl}).
To deduce the relation (\ref{eq:brelab1}) we use at first the sixth
relation in (\ref{eq:srebkl}):
\begin{equation*}
a_{(i+1)i} a_{(i+2)(i+1)} b_{(i+1)i} = a_{(i+1)i} b_{(i+2)i}
a_{(i+2)(i+1)},
\end{equation*}
and then the seventh relation in (\ref{eq:srebkl}):
\begin{equation*}
a_{(i+1)i} b_{(i+2)i} a_{(i+2)(i+1)} = b_{(i+2)(i+1)} a_{(i+1)i}
a_{(i+2)(i+1)}.
\end{equation*}
To deduce the relation (\ref{eq:brelab2}) we use the second, the
eighth, the fifth and again the second relations in
(\ref{eq:srebkl}):
\begin{multline*}
a_{(i+2)(i+1)} a_{(i+1)i} b_{(i+2)(i+1)} = a_{(i+1)i} a_{(i+2)i}
b_{(i+2)(i+1)}=
a_{(i+1)i} b_{(i+1)i} a_{(i+2)i} = \\
b_{(i+1)i} a_{(i+1)i}
a_{(i+2)i} = b_{(i+1)i} a_{(i+2)(i+1)}a_{(i+1)i}.
\end{multline*}
To eliminate (\ref{eq:defbab}) apply the second relation in
(\ref{eq:srebkl}) to change the center pair in 
(\ref{eq:defbab})
\break
$a_{(p+2)(p+1)} b_{(p+1)p}$ to
$b_{(p+2)p} a_{(p+1)p}$.
Then apply this process to the pair $a_{(p+3)(p+2)} b_{(p+2)p}$.
Ultimately, this process will move the original center letter
$b_{(p+1)p}$, to the leftmost position, where it becomes
$b_{qp}$. Free cancellation eliminates everything to its right,
and we are done.
\end{proof}
Now we consider the {\it positive singular braid monoid} with respect 
to generators $a_{ts}$ and $b_{t,s}$ for $1\leq s < t \leq n$. Its
relations are (\ref{eq:srebkl}) except one concerning the
invertibility of $a_{ts}$. Two positive words $A$ and $B$ in the
alphabet $a_{ts}$ and $b_{t,s}$
will be said to be {\it positively equivalent} if they are equal as 
elements of this monoid. In this case, as in the previous section, 
we shall write $A\doteq B$.

\begin{figure}
\epsfbox{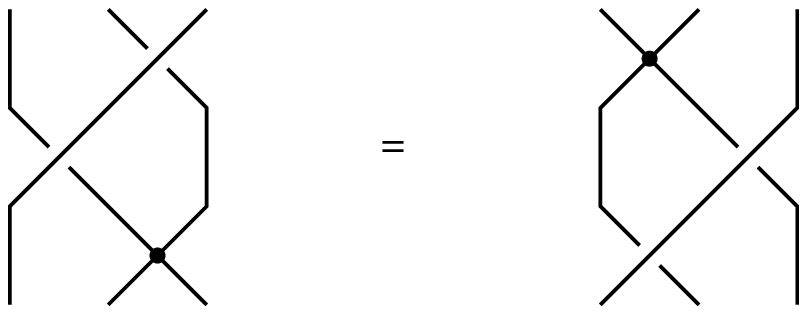}
\caption{}
\label{fi:bdelta}
\end{figure}

The {\it fundamental word} $\delta$ of Birman, Ko and Lee is given by
the formula
$$\delta \equiv a_{n(n-1)}a_{(n-1)(n-2)}
\dots a_{21} \equiv \sigma_{n-1} \sigma_{n-2} \dots \sigma_2\sigma_1.$$
Its divisibility by any generator $a_{ts}$, proved in \cite{BKL},
is convenient for us to be expressed in the following form.
\begin{pro} The fundamental word $\delta$ is positively equivalent to
 a word that begins or ends with any given generator $a_{ts}$.
The explicit expression for left divisibility is
\begin{equation*}
\delta \doteq a_{ts}a_{n(n-1)}a_{(n-1)(n-2)}
\dots a_{(t+1)s} a_{t(t-1)}\dots a_{(s+2)(s+1)} a_{s(s-1)}
\dots a_{21}.
\end{equation*}
\end{pro}
\begin{pro} For the fundamental word $\delta$ there are
the following formulae of commutation
\begin{equation*}
\begin{cases}
a_{ts} \delta &\doteq \delta a_{(t+1)(s+1)} \ \ \text{for}
\ \ 1\leq r < s < t < n, \\
a_{ns} \delta &\doteq \delta a_{(s+1)1}, \\
b_{ts} \delta &\doteq \delta b_{(t+1)(s+1)} \ \ \text{for} \ \
 1\leq r < s < t < n, \\
b_{ns} \delta &\doteq \delta b_{(s+1)1}.
\end{cases}
\end{equation*}
\end{pro}
\begin{proof} For the generators $a_{ts}$ it is proved in \cite{BKL}.
Suppose that $1\leq r < s < t < n$, then using
relations (\ref{eq:srebkl}) we have
\begin{multline*}
b_{ts} \delta \doteq b_{ts} a_{ts}a_{n(n-1)}a_{(n-1)(n-2)}
\dots a_{(t+1)s} a_{t(t-1)}\dots a_{(s+2)(s+1)} a_{s(s-1)}
\dots a_{21} \doteq \\
a_{ts} b_{ts}a_{n(n-1)}a_{(n-1)(n-2)} \dots a_{(t+1)s}
a_{t(t-1)}\dots a_{(s+2)(s+1)} a_{s(s-1)}\dots a_{21} \doteq \\
a_{ts} a_{n(n-1)}a_{(n-1)(n-2)} \dots b_{ts}a_{(t+1)s}
a_{t(t-1)}\dots a_{(s+2)(s+1)} a_{s(s-1)}\dots a_{21} \doteq \\
a_{ts} a_{n(n-1)}a_{(n-1)(n-2)} \dots a_{(t+1)s}b_{(t+1)t}
a_{t(t-1)}\dots a_{(s+2)(s+1)} a_{s(s-1)}\dots a_{21} \doteq \\
a_{ts} a_{n(n-1)}a_{(n-1)(n-2)} \dots a_{(t+1)s}
a_{t(t-1)}b_{(t+1)(t-1)}
\dots a_{(s+2)(s+1)} a_{s(s-1)}\dots a_{21} \doteq \\
a_{ts} a_{n(n-1)}a_{(n-1)(n-2)} \dots a_{(t+1)s}
a_{t(t-1)}\dots b_{(t+1)(s+2)}
a_{(s+2)(s+1)} a_{s(s-1)}\dots a_{21} \doteq \\
a_{ts} a_{n(n-1)}a_{(n-1)(n-2)} \dots a_{(t+1)s}
a_{t(t-1)}\dots a_{(s+2)(s+1)} b_{(t+1)(s+1)}
a_{s(s-1)}\dots a_{21} \doteq \\
a_{ts} a_{n(n-1)}a_{(n-1)(n-2)} \dots a_{(t+1)s}
a_{t(t-1)}\dots a_{(s+2)(s+1)}
a_{s(s-1)}\dots a_{21}b_{(t+1)(s+1)} \doteq \delta b_{(t+1)(s+1)}.
\end{multline*}
For $t= n$, we have
\begin{multline*}
b_{ns} \delta \doteq b_{ns} a_{n(n-1)} a_{(n-1)(n-2)} \dots a_{21}
\doteq  a_{n(n-1)}b_{(n-1)s} a_{(n-1)(n-2)} \dots a_{21} \doteq
 \dots \doteq \\
a_{(n-1)n} a_{(n-1)(n-2)}\dots a_{(s+2)(s+1)}b_{(s+1)s}
a_{(s+1)s} \dots a_{21} \doteq \\
a_{(n-1)n} a_{(n-1)(n-2)}\dots a_{(s+1)s}b_{(s+1)s}
a_{s(s-1)} \dots a_{21} \doteq  \\
a_{(n-1)n} a_{(n-1)(n-2)}\dots a_{s(s-1)}
b_{(s+1)(s-1)} a_{(s-1)(s-2)}\dots a_{21} \doteq  \\
a_{(n-1)n} a_{(n-1)(n-2)}\dots a_{s(s-1)}
 a_{(s-1)(s-2)}\dots a_{21}b_{(s+1)1}.
\end{multline*}
\end{proof}
Geometrically this commutation is shown on Figure~\ref{fi:bdelta}.

The analogues of the other results proved by Birman, Ko and Lee
(which are also analogues of the statements of Sections~2 and 3
of the present paper) remain valid for the singular braid monoid.
Unfortunately the proof of the analogue of Proposition~\ref{pro:div}
(which consists of many different cases already in \cite{BKL}) is very 
long and consists of even bigger amount of cases. 

\end{document}